\documentclass[final,12pt]{elsarticle}

\usepackage{a4wide}

\usepackage{hyperref}

\usepackage{amsmath}
\usepackage{amssymb}
\usepackage{amsthm}

\newtheorem{lemma}{Lemma}
\usepackage{xcolor}

\usepackage{multirow}

\usepackage{caption}
\usepackage{subcaption}

\usepackage{tikz}
\usepackage{circuitikz}
\usepackage{pgfplots}

\ctikzset{resistors/scale=0.8,
          capacitors/scale=0.7,
          diodes/scale=0.6,
          transistors/scale=1.3} 

\usepackage{hyperref}

\def\eps{\epsilon}
\def\ddt{\frac{d}{dt}}

\newcommand{\sign}{\text{sign}}

\begin{document}

\begin{frontmatter}

\title{MONA -- A magnetic oriented nodal analysis for electric circuits}

\author{Idoia Cortes$^a$ Garcia}\ead{idoia.cortes@tu-darmstadt.de}
\author{Herbert Egger$^b$}\ead{herbert.egger@ricam.oeaw.ac.at}
\author{Vsevolod Shashkov$^a$\corref{cor}}\ead{shashkov@mathematik.tu-darmstadt.de}

\cortext[cor]{Corresponding author}
\address{$^a$TU Darmstadt, Dolivostr. 15, 64293 Darmstadt, Germany}
\address{$^b$Johannes Kepler University Linz Johann Radon Institute for Computational and Applied Mathematics, Altenbergerstr. 69, A-4040 Linz, Austria}

\begin{abstract}
The modified nodal analysis (MNA) is probably the most widely used formulation for the modeling and simulation of  electric circuits. Its conventional form uses electric node potentials and currents across inductors and voltage sources as unknowns, thus taking an \emph{electric viewpoint}. In this paper, we propose a \emph{magnetic oriented nodal analysis } (MONA) for electric circuits, which is based on magnetic node potentials and charges across capacitors and voltage sources as the primary degrees of freedom, thus giving direct access to these quantities. 
The resulting system has the structure of a generalized gradient system which immediately ensures passivity in the absence of sources.
A complete index analysis is presented showing regularity of the magnetic oriented formulation under standard topological conditions on the network interconnection. In comparison to  conventional MNA, the differential-algebraic index is reduced by one in most cases which facilitates the numerical solution.   
Some preliminary numerical experiments are presented for illustration of the feasibility and stability of the new approach.
\end{abstract}


\end{frontmatter}


\section{Introduction}

The modeling and simulation of electric devices is one of the fundamental problems in electrical engineering. Very often the system dynamics can be described as an electric interconnection network \cite{campbell2019applications,gunther1995dae,gunther2005modelling}, by equivalent electric or magnetic circuit models \cite{amrhein20093,laithwaite1967magnetic,ruehli1974equivalent}, or by combinations of circuit and field equations \cite{bannenberg2021coupling,bartel2011structural,bartel2018pdaes,GarciaEtAl2020,tischendorf2003coupled}. 
Since its introduction in the mid 70s \cite{ho1975modified}, the \emph{modified nodal analysis} (MNA) has become the industry standard for electric circuit design and the simulation. We refer to \cite{gunther2005modelling,riaza2008differential} for an overview and further references.
In compact form, the governing equations read   
\begin{align}
    A_C \ddt q(A_C^\top e) + A_R g(A_R^\top e)  + A_L i_L + A_V i_V &= -A_I i_{src}, \label{eq:mna1} \\
 \ddt \phi(i_L) - A_L^\top  e &= 0, \label{eq:mna2} \\
 -A_V^\top  e &= -v_{src}. \label{eq:mna3}
\end{align}
Here $e$ denotes the vector of electric node potentials, $i_L$, $i_V$ are the currents through inductors and voltage sources, and $i_{src}$, $v_{src}$ are prescribed source terms. 
The interconnection of the individual devices is encoded in the partial incidence matrices $A_x$, with $x$ characterizing the type of the circuit element.
The characteristics of capacitors, inductors and resistors are described by the nonlinear device functions $q(v_C)$, $\phi(i_L)$, $g(v_R)$ with $v_x = A_x^\top  e$ denoting the voltages across elements. 
These nonlinear relations may be expressed as
\begin{align} \label{eq:mna4}
\ddt q(v_C) = C(v_c) \, \ddt v_c, \qquad 
\ddt \phi(i_L) = L(i_L) \, \ddt i_L, \qquad 
g(v_R) = G(v_R) \, v_R
\end{align}
with $C(v_C)$ and $L(i_L)$ denoting the more common differential capacitance and inductance matrices, respectively. 
These matrices are related to the Hessians of corresponding electric and magnetic energy functionals and, consequently, they can be assumed symmetric and positive definite; see Sections~\ref{sec:2} and \ref{sec:3} for details.
In a similar manner, the conductance matrix $G(v_R)$ characterizing resistive elements can be related to a dissipation functional and thus again assumed symmetric and positive definite. 

Substituting the identities \eqref{eq:mna4} into \eqref{eq:mna1}--\eqref{eq:mna3} leads to the \emph{conventional form} of the MNA. Its particular dissipative Hamiltonian structure automatically guarantees passivity of the system; see e.g. \cite{gunther2021dynamic,scholz2018signature} and Section~\ref{sec:2}.
A mathematical subtelty arising in the context of circuit modeling is the  differential-algebraic nature of the governing systems \cite{griepentrog1986,petzold1982differential}. 
In fact, the research in differential-algebraic equations and their numerical solution has been stimulated substantially for many years by applications in electronic circuits; see \cite{brenan1995numerical,hairer2006numerical,KunkelMehrmann2006} for an introduction and further references. 
It is now well understood by now that for consistent initial conditions and appropriate device characteristics
\begin{itemize}
\item[(a1)] the system \eqref{eq:mna1}--\eqref{eq:mna4} is well-posed and has index $\nu \le 2$ if {the circuit contains neither loops of voltage sources nor cutsets of current sources}.
\end{itemize}
Moreover, one can show that
\begin{itemize}
\item[(a2)] the index is $\nu \le 1$ if {the circuit contains neither loops of capacitors and voltage sources nor cutsets of inductors and current sources};
\end{itemize}
see \cite{estevez2000structural,gunther1995dae, tischendorf1999topological} for details and proofs. 
The condition in (a2) could in fact be slightly relaxed; cf. \cite[Remark~***]{gunther2021dynamic}.
Let us note that the conditions in (a1) and (a2) are purely topological, i.e., concerning only the interconnection of devices, and they can therefore be formulated as algebraic conditions on the partial incidence matrices $A_x$; see \cite{estevez2000structural,tischendorf1999topological} for an explanation. 
This allows for a systematic \emph{projection-based analysis} \cite{griepentrog1986,lamour2013} and different index concepts have been successfully employed to prove (a1)--(a2) and generalizations; we refer to \cite{lamour2013,reis2014mathematical,riaza2008differential} for further results and references.

Besides their analytical peculiarities, electric circuit equations also pose various challenges for the numerical solution. 
Due to the differential-algebraic nature, implicit time stepping schemes have to be used \cite{brenan1995numerical,hairer2006numerical,KunkelMehrmann2006}.
While passivity on the discrete level can be proven rigorously for the implicit Euler method and related variational time discretization schemes of higher order \cite{egger2019structure}, 
strict passivity may in general be lost through discretization by standard single or multistep schemes. 
In the presence of strong nonlinearities, even well-established second order schemes, like the Crank-Nicolson or BDF-2 method may run into stability problems.
A common practice in industry therefore is to use low order time integration schemes in general and eventually fall back to the implicit Euler method in case of stability issues. 
Let us refer to \cite[Ch.~10,11]{gunther2005modelling}
for details. 

Another somehow related difficulty is that the nonlinear differential relations \eqref{eq:mna4} for device characteristics, which were employed in the derivation of the conventional form of the MNA, can in general not be reproduced exactly on the discrete level. 
As a consequence, charge conservation may be lost after discretization and simulated magnetic fluxes may be inconsistent. 
A possible remedy is to introduce extra variables 
\begin{align} \label{eq:mna5}
q_C = q(v_C), \qquad \phi_L = \phi(i_L)
\end{align}
for the electric charges and magnetic fluxes which together with \eqref{eq:mna1}--\eqref{eq:mna3} leads to the \emph{charge/flux oriented MNA}; see \cite{gunther2005modelling} for an overview.
The assertions (a1) and (a2) can be verified also for this extended formulation \cite{estevez2000structural} but some additional modifications are required to regain the dissipative Hamiltonian structure \cite{gunther2021dynamic}. 
Such extensions typically involve substantially more unknowns than the conventional MNA, in particular, if many capacitative or inductive devices are present. 

Let us emphasize that all formulations mentioned so far involve the electric node potentials $e$ as the primary unkowns and could therefore be called \emph{electric oriented nodal analyses}.
In this paper, we take an alternative \emph{magnetic viewpoint} and introduce a \emph{magnetic oriented nodal analysis} (MONA) for the modeling and simulation of electric circuits. 
In its compact form, the resulting formulation reads
\begin{align}
    A_R g(A_R^\top  \ddt \psi)  + A_C \ddt q_C + A_V \ddt q_V  
    + A_L i(A_L^\top \psi) 
    &=- A_I i_{src}, \label{eq:mona1}\\
    -A_C^\top  \ddt \psi 
    + v(q_C) &= 0, \label{eq:mona2}\\
    -A_V^\top  \ddt \psi &= -v_{src}. \label{eq:mona3}
\end{align}
Here $\psi$ denotes a vector of magnetic node potentials and $q_C,q_V$ are the displaced charges at capacitors and voltage sources. 
The formulation obviously has a great similarity with the MNA, allowing to reuse available implementations, but the new formulation also has some subtle differences.
Access to the magnetic flux linkages and electric node potentials is now available through 
\begin{align} \label{eq:mona4}
\phi_L = A_L \psi, \qquad e = \ddt \psi,
\end{align}
from which voltages and currents can be derived like in the electric based nodal analysis above. 
While the conductance relation $g(v_C) = G(v_C) \, v_C$ is the same as in \eqref{eq:mna1}, the device characteristics of inductors and capacitors are now described by 
\begin{align} \label{eq:mona5}
i(\phi_L) = \nabla \eps_L(\phi_L), \qquad 
v(q_C) = \nabla \eps_C(\phi_C)
\end{align}
with $\eps_L(\cdot)$, $\eps_C(\cdot)$ denoting magnetic and electric energy functionals, respectively. 
To draw the connection with the previous models, we may differentiate equation \eqref{eq:mona5} to see that
\begin{align} \label{eq:mona6}
    \ddt i(\phi_L) = \nabla^2 \eps_L(\phi_L) \, \ddt \phi_L, 
    \qquad 
    \ddt v(q_C) = \nabla^2 \eps_C(q_C) \, \ddt q_C.
\end{align}
A quick comparison with the relations in \eqref{eq:mna4} reveals that the Hessians $\nabla^2 \eps_L(\phi_L) = L(i(\phi_L))^{-1}$ and $\nabla^2 \eps_C(q_C) = C(v(q_C))^{-1}$ are directly linked to the differential inductance and capacitance matrices employed before.
The derivation of the system \eqref{eq:mna1}--\eqref{eq:mna3} is based on the very same principles as the new model \eqref{eq:mona1}--\eqref{eq:mona3}, and the electric and magnetic viewpoints are therefore mutually equivalent and allow to model the same circuits.

Let us note that in the magnetic oriented formulation \eqref{eq:mona1}--\eqref{eq:mona4}, the dissipative terms and the interconnection structure affect the time derivatives, in contrast to \eqref{eq:mna1}--\eqref{eq:mna3}, where they involve only the lowest order terms. 
The system \eqref{eq:mona1}--\eqref{eq:mona4} therefore has a different geometric structure as the MNA, namely that of a \emph{generalized gradient system}, which however again immediately allows to guarantee passivity of the system; see  \cite{egger2020energy,hairer2014} and Section~\ref{sec:3} below.
We will further show that for consistent initial conditions and appropriate device characteristics
\begin{itemize}
    \item[(b1)] the system \eqref{eq:mona1}--\eqref{eq:mona3} is well-posed and has index $\nu \le 1$ if the circuit contains neither loops of voltage sources nor cutsets of current sources;
    \item[(b2)] the index is $\nu=0$ if the circuit contains neither loops of capacitors and voltage sources nor cutsets of inductors and current sources.
\end{itemize}
Let us emphasize that the conditions in (b1) and (b2) are 
exactly the same
as those employed in the assertions (a1) and (a2) before, i.e., the new magnetic oriented formulation is as flexible and general as the electric one. 
In comparison to the conventional MNA, the proposed magnetic oriented formulation however leads to an index reduction in most cases, which alleviates the numerical solution to some extent. 
Apart from this numerical advantage, the main benefit of our our approach certainly is the direct access to the charges $q_C$, $q_V$ and flux linkages $\phi_L=A_L^\top \psi$, which might be of particular interest in certain applications, e.g., if many energy storing elements are present in the circuit.
A disadvantage of the magnetic oriented formulation, on the other hand, seems to be that access to electric quantities, like electric potentials and branch currents, is somewhat indirect, e.g., via $e = \ddt \psi$, which however can be realized exactly also on the discrete level.

The remainder of the manuscript is organized as follows:
In Section~\ref{sec:2}, we introduce our notation and recall some basic facts about the conventional electric oriented nodal analysis.
The new magnetic oriented formulation is then derived in Section~\ref{sec:3} and a short proof of its passivity is provided. 
Section~\ref{sec:4} is concerned with the index analysis of MONA and contains the proof of assertions (b1) and (b2). 
Some preliminary numerical results for our method and a brief comparison with the conventional MNA are presented in Section~\ref{sec:5}.

\section{Basic notation and a review of MNA} \label{sec:2}

We consider a directed and connected graph with $N_n$ nodes and $N_b$ branches. Its interconnection structure is described by the \textit{incidence matrix} $A\in \mathbb{R}^{N_n\times N_b}$ defined as
\begin{align*}
A_{ij} = 
\begin{cases}
\phantom{-}1 & \text{if branch } j \text{ leaves node } i,\\
-1& \text{if branch } j \text{ enters node } i, \\
\phantom{-}0 & \text{else}.
\end{cases}
\end{align*}
We use the same letter $A$ to denote the reduced incidence matrix which results from eliminating one row corresponding to a grounded vertex. 
Kirchhoff's current law and the definition of the voltages across elements can then be expressed as 
\begin{equation}
A i = 0 \quad\text{and} \quad v = A^\top e,
\label{eq:Kirchhoffs}
\end{equation}
with $i$, $v$, and $e$ denoting the vectors of branch currents, branch voltages and node potentials, respectively.
We consider circuits consisting of resistors, inductors, capacitors, and independent current and voltage sources. This allow to split 
$$
A = [A_C,A_L,A_R,A_V,A_I]
$$  
with $A_x$ representing the partial incidence matrix of the elements of type $x$.
In a similar manner, we split the vectors of currents and voltages and denote by $i_x$ and $v_x$ the corresponding sub-vectors. 
The system \eqref{eq:Kirchhoffs} is complemented by consitutive relations
\begin{align}
i_C =  C(v_C) \ddt v_C, \qquad 
v_L = L(i_L) \ddt i_L,  \qquad \text{and}\qquad 
i_R = G(v_R) v_R,
\label{eq:material_laws}
\end{align}
which describe the device characteristics of capacitors, inductors, and resistors, respectively. 
The matrices $C(v_C)$, $L(i_L)$, and $G(v_R)$ are the differential capacitance, differential inductance, and the generalized conductance, respectively, and they are assumed to be symmetric positive definite and to depend smoothly on their arguments.  
The currents $i_I=i_{src}$ and voltages $v_V=v_{src}$ denote the input to the system. 
Substituting \eqref{eq:material_laws} into \eqref{eq:Kirchhoffs} leads to the system
\begin{align}
A_C C(A_C^\top e)A_C^\top  \ddt e + A_R G(A_R^\top e)A_R^\top  e + A_L i_L + A_V i_V &= -A_I i_{src} \label{eq:eona1}\\
L(i_L) \ddt i_L - A_L^\top  e &= 0 \label{eq:eona2}\\
-A_V^\top  e &= -v_{src} \label{eq:eona3}
\end{align}
which is the conventional form of the modified nodal analysis \cite{ho1975modified}; also see \cite{estevez2000structural,gunther2005modelling}. 
In the following, we briefly recall the most important results of its analysis.

\bigskip 

\noindent
\textbf{Energy balance and passivity.} 
We start with deriving the basic energy--dissipation identity. 
To this end, let $\eps_L(\phi_L)$ and $\eps_C(q_C)$ denote the energy stored in inductors and capacitors, respectively, and recall \cite{estevez2000structural,gunther2005modelling} that
\begin{align*}
    \nabla \eps_L(\phi_L) = i_L, \qquad \nabla \eps_C(q_C) = v_C.
\end{align*}
With these identities and expressing 
$\phi_L=\phi(i_L)$ and $q_C=q(v_C)$, we immediately obtain 
\begin{align*}
    \ddt \Big(  \eps_C(q(v_C)) + \eps_L(\phi(i_L)) \Big) 
    &=\langle \nabla \eps_C(q(v_C)), \ddt q(v_C) \rangle +  \langle \nabla \eps_L(\phi(i_L)), \ddt \phi(i_L) \rangle\\
    &= \langle v_C, C(v_C) \, \ddt v_C\rangle + \langle i_L, L(i_L) \, \ddt i_L \rangle, 
\end{align*}
where we used \eqref{eq:mna4} in the second step.
Here and below, $\langle a, b \rangle = b^\top a$ is the Euclidean inner product.
After substituting $v_C=A_V^\top e$, one can see that the result amounts to the first two terms in \eqref{eq:eona1}--\eqref{eq:eona2} multiplied by $e^\top$ and $i_L^\top$ from the left. 
Inserting these equations and rearranging the terms thus leads to 
\begin{align*}
  \ddt \left(  \eps_C(q(v_C)) + \eps_L(\phi(i_L)) \right) 
  &= -\langle e, A_R^\top G(v_R) A_R^\top e + A_L i_L + A_V i_V + A_I i_{src} \rangle + \langle i_L, A_L^\top e\rangle \\ 
  &= -\langle e, A_R G(v_R) A_R^\top e \rangle  - \langle A_I^\top e, i_{src}\rangle - \langle i_V, v_{src}\rangle,
\end{align*}
where we used \eqref{eq:eona3} and some elementary algebraic manipulations in the last step. 
This identity states that the energy of the system changes only by dissipation in resistors and supply or loss through voltage and current sources.   
In particular, when $i_{src}=v_{src}=0$, we obtain passivity of the system.

\bigskip 

\noindent
\textbf{Index analysis.} 
As a second step, let us recall some basic facts about the index analysis of the MNA equations. Recall that $L(i_L)$, $C(v_C)$ and $G(v_R)$ are assumed symmetric and positive definite matrices that depend smoothly on their arguments. 
\begin{lemma} \label{lem:mna}
Assume that 
\begin{align} \label{eq:A1}
    N([A_R,A_C,A_V,A_L]^\top)=0 \qquad  \text{and} \qquad N(A_V)=0.  \tag{A1}
\end{align}
Then \eqref{eq:eona1}--\eqref{eq:eona3} is a regular system of DAEs with index $\nu \le 2$. 
If additionally
\begin{align} \label{eq:A2}
    N([A_R,A_C,A_V]^\top)=0 \qquad \text{and} \qquad N([A_C,A_V])=0. \tag{A2} 
\end{align}
Then the system is again regular and of index $\nu \le 1$. 
\end{lemma}
\noindent
A detailed proof of these results can be found in various references; see e.g. \cite{estevez2000structural,reis2014mathematical,tischendorf1999topological}. 
Let us note that the algebraic conditions (A1)--(A2) are equivalent to the topological conditions in (a1)--(a2) mentioned in the introduction; see again \cite{estevez2000structural,gunther2005modelling} for details and explanation. 
Like in assertion (a2) the condition (A2) could be slightly relaxed for the second statement.

\section{The magnetic oriented nodal analysis}\label{sec:3}

By Faraday's law, the voltage induced by a time varying magnetic flux through a wire loop is given by $v_L = \ddt \phi_L$. Integrating this expression and using the definition $v_L=A_L^\top e$ of the voltages across inductors, we obtain $A_L^\top \int_0^t e dt = \phi_L$, where we assumed $e$ to be zero at time $t=0$. 
We now introduce a \emph{magnetic node potential} by 
$\psi = \int_0^t e \, dt$, which gives 
\begin{align}
    e = \ddt \psi, \qquad \phi_L = A_L^\top \psi.    
\end{align}
In a similar spirit, we consider an integral form of charge conservation by introducing generalized charges $q_x = \int_0^t i_x \, dt$. This allows us to compute the currents and voltages
\begin{align}
i_C = \ddt q_C, \quad i_V = \ddt q_V 
\qquad \text{and} \qquad v_C = A_C^\top \ddt \psi, \quad v_V = A_V^\top \ddt \psi
\end{align}
across capacitors and voltage sources, respectively. 
Inserting these expressions in Kirchhoff's law \eqref{eq:Kirchhoffs} and adding equations $v_C=v(q_C)$ and $v_V=v_{src}$ for the voltages across capacitors and voltage sources now immediately leads to the equations \eqref{eq:mona1}--\eqref{eq:mona3}. 
Employing the relations $g(v_R)=G(v_R) \, v_R$ and \eqref{eq:mona5} for the device characteristics, we arrive at 
\begin{align}
    A_R G(A_R^\top \ddt \psi) A_R^\top \ddt \psi + A_C \ddt q_C + A_V \ddt q_V &= -A_L \nabla \eps_L(A_L^\top \psi) - A_I i_{src} \label{eq:mona1b}\\
    -A_C^\top \ddt \psi &= -\nabla \eps_C(q_C) \label{eq:mona2b}\\
    -A_V^\top \ddt \psi &= -v_{src} \label{eq:mona3b},
\end{align}
which can be considered as the magnetic oriented formulation of \eqref{eq:eona1}--\eqref{eq:eona3}.
Let us note that the electric and magnetic oriented formulations are derived on the basis of the same physical principles and they can be transformed into each other.

\bigskip 
\noindent\textbf{Energy balance and passivity.}
As a first step of our analysis, we derive a basic energy-dissipation identity. 
By formal differentiation and noting that $\phi_L = A_L^\top \psi$, we see that
\begin{align*}
    \ddt \left( \eps_L(\phi_L) + \eps_C(q_C)\right)
    &= \langle \nabla \eps_L(\phi_L), A_L^\top \ddt \psi\rangle + \langle \nabla \eps_C(q_C), \ddt q_C\rangle.
\end{align*}
Apart from the source current, this amounts to the right hand sides of \eqref{eq:mona1b}--\eqref{eq:mona2b} multiplied by $\ddt \psi^\top$ and $\ddt q_C^\top$ from the left. Inserting the equations therefore gives 
\begin{align*}
    &\ddt \left( \eps_L(\phi_L) + \eps_C(q_C)\right)
    = -\langle A_R G(A_R^\top \ddt \psi) A_R^\top \ddt \psi + A_C \ddt q_C + A_V \ddt q_V + A_I i_{src},  \ddt \psi\rangle  \\
    &\qquad\qquad\qquad + \langle A_C^\top \ddt \psi, \ddt q_C\rangle 
    = -\langle A_R G(A_R^\top \ddt \psi) A_R^\top \ddt \psi, \ddt\psi \rangle - \langle i_{src}, v_I\rangle - \langle i_V, v_{src}\rangle, 
\end{align*}
where we used \eqref{eq:mona3b} and some elementary algebraic manipulations in the last step. 
The energy of the system can therefore again only change due to dissipation in the resistors and power supplied or drawn through current and voltage sources. 
In particular, the system is passive in the absence of source terms.

\section{Index analysis}\label{sec:4}

As a second step of our analysis, we now consider in detail the differential-algebraic index of the system under consideration. For ease of notation,  we choose $i_{src}=v_{src}=0$ and then rewrite the system \eqref{eq:mona1b}--\eqref{eq:mona3b} in compact form 
\begin{align} \label{eq:mona}
    \begin{pmatrix}
    A_R G^\psi A_R^\top  & A_C & A_V \\
    -A_C^\top  \\  -A_V^\top  \end{pmatrix}
    \ddt 
    \begin{pmatrix} \psi \\ q_C \\ q_V \end{pmatrix}
    = - \begin{pmatrix}
      A_L \nabla \eps_L(A_L^\top  \psi) \\
      \nabla \eps_C(q_C) \\ 0 \end{pmatrix}
\end{align}
We use $G^\psi=G(A_R^\top \ddt \psi)$ to abbreviate the solution dependent conductivity matrix in the following. 
In addition, we assume regular device characteristics, which can be stated as 
\begin{align} \label{eq:A0} \tag{A0}
    \text{$G^\psi$ symmetric and positive definite and $\eps_L$, $\eps_C$ are strictly convex,} 
\end{align}
with all functions depending smoothly on their arguments. 
Let us note that these were the standard assumptions also for the investigation of the MNA summarized in the previous section. 
Under the algebraic conditions of Lemma~\ref{lem:mna}, we then obtain the following result.

\bigskip 

\begin{lemma} \label{lem:mona}
Let \eqref{eq:A0} hold and assume that 
\eqref{eq:A1} is valid, i.e., 
\begin{align*}
    N([A_R,A_C,A_V,A_L]^\top)=0 \qquad \text{and} \qquad N(A_V)=0.
\end{align*}
Then the system \eqref{eq:mona} is a regular DAE of index $\nu=1$.
If \eqref{eq:A2} holds additionally to \eqref{eq:A1}, i.e.,
\begin{align*}
N([A_R,A_C,A_V]^\top)=0 \qquad \text{and} \qquad N([A_C,A_V])=0,
\end{align*}
then \eqref{eq:mona} is a regular DAE of index $\nu=0$, i.e., an ordinary differential equation.
\end{lemma}
\begin{proof}
We start with the second assertion: 
Define $B=[A_C,A_V]^\top$, which is surjective since 
$$
R(B)^\perp 
= R([A_C,A_V]^\top)^\perp =N([A_C,A_V]) = 0.
$$ 
Further set $C=A_R G^\psi A_R^\top$ and note that
for $\psi \ne 0$ we have 
$$
\langle C \psi, \psi\rangle 
= \langle G^\psi A_R^\top \psi, A_R^\top \psi\rangle 
 \ge \gamma \|A_R^\top \psi\|^2 \ne 0, 
$$ 
with some $\gamma>0$, since $G^\psi$ is positive definite by assumption \eqref{eq:A0}. 
By assumption \eqref{eq:A2}, we further see that $N(A_R^\top) \cap N([A_C,A_V]^\top) = N([A_R,A_C,A_V]^\top) = 0$ and hence
\begin{align*}
\langle C \psi, \psi\rangle \ne 0 \qquad \text{for all } \psi \in N(B), \, \psi \ne 0.
\end{align*}
By Brezzi's splitting lemma \cite{brezzi1974existence}, we conclude that the leading matrix in \eqref{eq:mona} is regular and hence the system is an implicit ODE, i.e., a regular DAE with index $\nu=0$.

We now turn to the first assertion of the lemma,
for which we use a typical projection-based analysis \cite{lamour2013}.
If the conditions \eqref{eq:A2} are not valid, we simply split the vector spaces
\begin{align*}
    V_\psi &= N([A_R,A_C,A_V]^\top ) \oplus N([A_R,A_C,A_V]^\top )^\perp  \\
    V_q    &= N([A_C,A_V]) \oplus N([A_C,A_V])^\perp
\end{align*}
of magnetic potentials $\psi$ and charges $q=(q_C,q_V)$.
By choosing orthogonal bases for the corresponding subspaces, we can thus decompose 
\begin{align*}
    \psi &= Q_1 \psi_1 + Q_2 \psi_2 
    \qquad \text{and} \quad q = P_1 q_1 + P_2 q_2.
\end{align*}
with stacked matrices $Q=[Q_1,Q_2]$ and $P=[P_1,P_2]$ being square and orthogonal. 
The second splitting can also be written component-wise as
\begin{align*}
    \begin{pmatrix} q_C \\ q_V \end{pmatrix}
    = 
    \begin{pmatrix} P_{1,C} \\ P_{1,V} \end{pmatrix} p_1 + 
    \begin{pmatrix} P_{2,C} \\ P_{2,V} \end{pmatrix} p_2.
\end{align*}
We now multiply the system \eqref{eq:mona} from left by $\text{blkdiag}(Q^\top,P^\top)$ and use the above expansions for $\psi$ and $q=(q_c,q_v)$,
which leads to the equvivalent form 
\begin{align} \label{eq:monab}
    \begin{pmatrix}
    0 & 0             & 0 & 0 \\
    0 & \tilde G_{22} & 0 & \tilde B\\
    0 & 0             & 0 & 0 \\
    0 & -B            & 0 & 0  
    \end{pmatrix}
    \ddt 
    \begin{pmatrix} \psi_1 \\ \psi_2 \\ q_1 \\q_2 \end{pmatrix}
    = - 
    \begin{pmatrix}
    Q_1^\top  A_L \nabla \eps_L(A_L^\top  (Q_1 \psi_1 + Q_2 \psi_2)) \\
    Q_2^\top  A_L \nabla \eps_L(A_L^\top  Q_1 \psi_1 + Q_2 \psi_2)) \\
    P_{1,C}^\top  \nabla \eps_C(P_{1,C} q_1 + P_{2,C} q_2) \\
    P_{2,C}^\top  \nabla \eps_C(P_{1,C} q_1 + P_{2,C} q_2) 
    \end{pmatrix}
\end{align}
By a simple rearrangement of the variables and the equations, this problem can then be written compactly as a Hessenberg system of the form
\begin{align*}
    M^\psi \ddt y &= f(y,z) \\
                0 &= g(y,z)
\end{align*}
with differential and algebraic variables $y=(\psi_2,q_2)$ and $z=(\psi_1,q_1)$, respectively.
Since all transformations employed so far are equivalence transformations, the rank of the matrix $M^\psi$ in front of the time derivative $\ddt y$ is the same as that of the leading matrix in \eqref{eq:mona}. 
A simple counting argument then implies that $M^\psi$ is regular. 
Let us now have a closer look onto the algebraic constraint 
\begin{align*}
    0 = g(y,z) := \begin{pmatrix} 
    Q_1^\top  A_L \nabla \eps_L(A_L^\top  (Q_1 \psi_1 + Q_2 \psi_2)) \\
    P_{1,C}^\top  \nabla \eps_C(P_{1,C} q_1 + P_{2,C} q_2) 
\end{pmatrix}
\end{align*}
By the chain rule, we obtain 
\begin{align} \label{eq:Dzg}
    D_z g(y,z) = 
    \begin{pmatrix} Q_1^\top  A_L \nabla^2 \eps_L(v_L) A_L^\top  Q_1 & 0 \\
    0 & P_{1,C}^\top  \nabla^2 \eps_C(q_C) P_{1,C} \end{pmatrix} 
\end{align}
with $v_L=A_L^\top (Q_1 \psi_1 + Q_2 \psi_2)$ and $q_C=P_{1,C} q_1 + P_{2,C} q_2$. 
In \eqref{eq:A0} the energy functionals are assumed strictly convex, and their Hessians $L^{-1} = \nabla^2 \eps_L(v_L)$ and $C^{-1} = \nabla^2 \eps_C(q_C)$ are thus symmetric and positive definite. 
The condition $N(A_L^\top ) \cap N([A_R,A_C,A_V]^\top ) = 0$ then implies that $A_L^\top  Q_1$ is injective and hence the upper left matrix in \eqref{eq:Dzg} is symmetric and positive definite. 
From the condition $N(A_V) = 0$, on the other hand, we can further deduce that $(0,q_V) \in N([A_C,A_V])$ implies $q_V=0$. 
This shows that the matrix $P_{1,C}$ is injective and hence the lower right matrix is again symmetric and positive definite. 
Hence the Jacobian $D_z g(y,z)$ is regular.
As a consequence, the above Hessenberg system is a regular DAE of index $\nu=1$. 
Since this system was obtained from  \eqref{eq:mona} by algebraic equivalence transformation, the same is true for the original system.
\end{proof}

\bigskip 
\noindent
\textbf{Summary.}
Let us close this section by briefly summarizing our theoretical observations: Under standard assumptions, the proposed magnetic oriented nodal analysis (MONA) is passive and a regular system of differential algebraic equations. In comparison to 
conventional MNA or the charge-flux based formulations of \cite{gunther2021dynamic}, the index is reduced by one in most cases. 

\section{Numerical illustration} \label{sec:5}

We now present numerical results for two simple test problems which demonstrate the stability and performance of MONA and highlight some possible advantages in comparision with the convential modified nodal analysis. 

\bigskip 

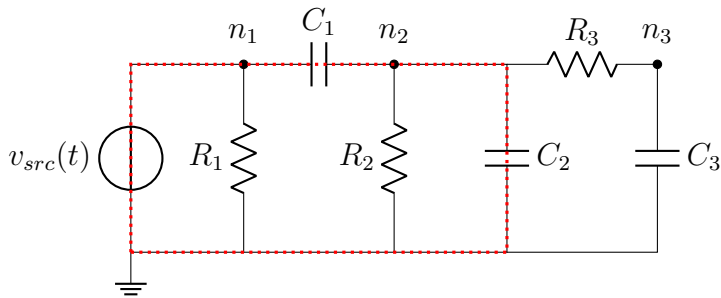
\begin{figure}[ht!]
    \centering
    \begin{circuitikz}
    \draw (0,0) to [V, l^=$v_{src}(t)$, -] (0,2.5) to (1.5,2.5) node[label=$n_1$]{} to [R, l_=$R_1$,*-] (1.5,0) to (0,0);
    \draw (1.5,2.5) to [C, l^=$C_1$, -] (3.5,2.5) to [R, l_=$R_2$,-] (3.5,0) to (1.5,0);
    \draw (3.5,2.5) node[label=$n_2$]{}  to [short,*-] (5,2.5)  to [C,l^=$C_2$,-] (5,0) to (3.5,0);
    \draw (5,2.5) to [R,l^=$R_3$,-*] (7,2.5) node[label=$n_3$]{} to [C, l=$C_3$, -] (7,0) to (5,0);
    \draw (0,0) to (0,-0.0) node[ground]{}; \draw[color=red,dotted,very thick] (0,0) to (0,2.5) to (5,2.5) to (5,0) to (0,0);
    \end{circuitikz}
    \caption{Index-2 circuit example from \cite[Fig.~10.2]{gunther2005modelling} containing a CV-loop marked with red dots.}
    \label{fig:circuit}
\end{figure}

\noindent\textbf{Example~1: An Index-2 problem.}
In our first example, which is taken from \cite[Ch.~10]{gunther2005modelling}, we consider an electric circuit containing a CV-loop; see Figure~\ref{fig:circuit} for an illustration. 
The conventional MNA therefore leads to an index-2 problem, while according to Theorem~\ref{lem:mona} the MONA approach results in a system with index 1. 
The circuit contains only capacitors, resistors and one voltage source, and the network topology is described by following partial incidence matrices
$$
A_C = 
\begin{pmatrix}
1 & 0 & 0 \\ -1 & 1 & 0 \\ 0 & 0 & 1
\end{pmatrix}
,\quad A_R = 
\begin{pmatrix}
1 & 0 & 0 \\ 0 & 1 & 1 \\ 0 & 0 & -1
\end{pmatrix} 
\quad
\mbox{and}
\quad
A_V = 
\begin{pmatrix}
1 \\ 0 \\ 0
\end{pmatrix}.
$$
For ease of presentation, we consider linear constitutive equations with $C_i=1$ and $R_i=1$ for $i=1,2,3$ and denote by $C=G=I_3$ the identity matrices of dimension $3$.
The conventional MNA formulation then leads to the following system 
$$
\begin{pmatrix}
A_C C A_C^T & 0 \\ 0 & 0
\end{pmatrix}
\ddt
\begin{pmatrix}
e \\ i_V
\end{pmatrix}
+
\begin{pmatrix}
A_R G A_R^T & A_V \\ -A_V^T & 0
\end{pmatrix}
\begin{pmatrix}
e \\ i_V
\end{pmatrix}
=
\begin{pmatrix}
0\\-v_{src}(t)
\end{pmatrix}
$$
while the magnetic-oriented scheme obtained by MONA reads
$$
\begin{pmatrix}
A_R G A_R^T & A_C & A_V \\ -A_C^T & 0 & 0 \\ -A_V^T &0 &0
\end{pmatrix}
\ddt
\begin{pmatrix}
\psi \\q_C \\q_V
\end{pmatrix}
+
\begin{pmatrix}
0 & 0 & 0 \\ 0 & C^{-1} & 0 \\ 0 & 0 & 0 
\end{pmatrix}
\begin{pmatrix}
\psi \\ q_C \\ q_V
\end{pmatrix}
=
\begin{pmatrix}
0 \\0\\-v_{src}(t)
\end{pmatrix}.
$$
Due to our choice of the constitutive equations, both systems are linear and time-invariant. 

\begin{figure}[ht!]
\centering
    \begin{tikzpicture}[scale = 0.75]
    \begin{axis}[
      grid = both,
      legend pos=north east
    ]
    \addplot[color = red,thick] table[x = t, y = ee1] {sine.dat};
    \addplot[color = green!50!black,thick] table[x = t, y = ee2] {sine.dat};
    \addplot[color = blue!90!black,thick] table[x = t, y = ee3] {sine.dat};
    \addplot[dotted,very thick] table[x = t, y = e1] {sine.dat};
    \addplot[dotted,very thick] table[x = t, y = e2] {sine.dat};
    \addplot[dotted,very thick] table[x = t, y = e3] {sine.dat};
    \legend{$e_1$,$e_2$,$e_3$};
    \end{axis}
    \end{tikzpicture}
    \hspace*{3em}
    \begin{tikzpicture}[scale = 0.75]
    \begin{axis}[
      grid = both,
      legend pos=north east
    ]
    \addplot[color = blue] table[x = t, y = iv] {sine.dat};
    \addplot[color = red,style={very thick}] table[x = t, y = iiv] {sine.dat};
    \legend{$i_V{MNA}$,$i_V^{MONA}$};
    \end{axis}
    \end{tikzpicture}
    \caption{Numerical solutions obtained by the trapezoidal rule applied to the MNA and MONA formulations. Left: potentials (MNA: dotted; MONA: solid); right: current through voltage source.}
    \label{fig:instability}
\end{figure}
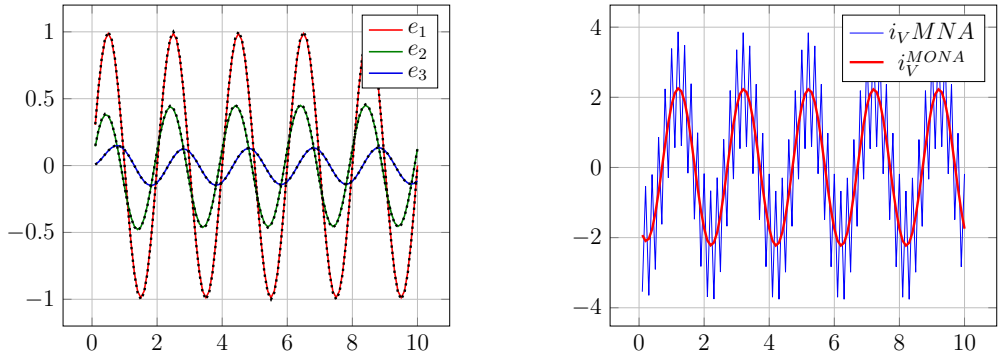

In Figure~\ref{fig:instability}, we display the electric solution components obtained by numerical solution of the two equations by the trapezoidal rule (TR) with a fixed time step $\tau = 0.1$ and for $v_{src}(t) = \sin(\pi t)$. 
For our simulations, we chose trivial initial conditions, which are consistent with the algebraic constraint caused by the voltage source. This suffices to guarantee stability for the index-1 formulation obtained by MONA. The MNA system, on the other hand, has index 2 and a \emph{hidden constraint} arises, which is not satisfied by our choice of initial conditions and causes large oscillations in the algebraic solution component. 
Let us note that this weak instability could be cured by an appropriate initialization phase, e.g., by performing the first time step with the implicit Euler method. 
If we choose the source term $v_{src}(t)=\cos(t)$ inconsistent with the trivial initial values, then the TR-discretization of the MNA formulation leads to strong instabilities, which can no longer be cured by one or two implicit Euler steps in the start-up phase of the simulation. 
In this case, also MONA shows a weak instability that can however be overcome by an appropriate start-up.

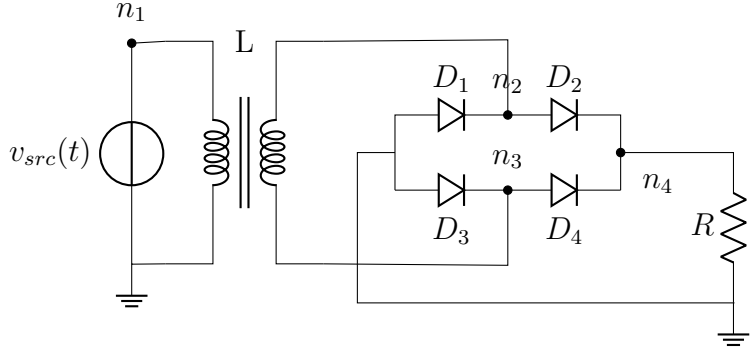
\begin{figure}[ht!]
    \centering
    \begin{circuitikz}
    \draw (1.5,1) node [transformer core, yscale=1.41](T){}
      (T.A1) node[above] {}
      (T.A2) node[below] {}
      (T.B1) node[above] {} 
      (T.B2) node[below] {}
      (T.base) node{L}; 
    \draw (T.A1) -- (0,2.455) node[label=$n_1$]{} to [V,l_=$v_{src}(t)$,*-] (0,-0.485) -- (T.A2);
    \draw (3.5,0.5) to [diode, l_=$D_3$] (5,0.5) to [diode, l_=$D_4$] (6.5,0.5) ;
    \draw (3.5,0.5) -- (3.5,1.5) to [diode, l^=$D_1$] (5,1.5) to [diode, l^=$D_2$] (6.5, 1.5) -- (6.5,0.5);
    \draw (T.B2) -- (5,-0.5)  to[short,-*] (5,0.5) node[label=above:$n_3$]{};
    \draw (T.B1) -- (5,2.5) to[short,-*] (5,1.5) node[label=north:$n_2$]{};
    \draw (3.5,1) -- (3,1) -- (3,-1) -- (8,-1);
    \draw (6.5,1)node[label=south east:$n_4$]{} to[short,*-] (8,1) to [R, l_=$R$] (8,-1);
    \draw (0,-0.485) to (0,-0.485) node[ground]{};
    \draw (8,-1) to (8,-1) node[ground]{};
    \end{circuitikz}
    \caption{Schematic sketch of a full wave rectifier.}
    \label{fig:rectifier}
\end{figure}

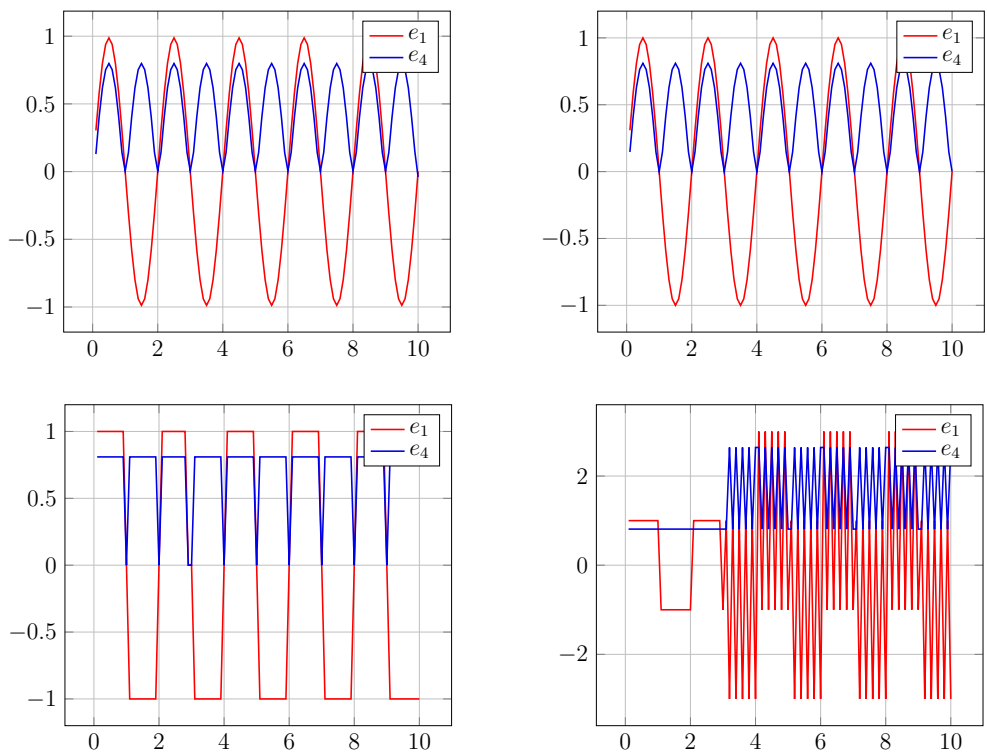
\begin{figure}[ht!]
\centering
   \begin{tikzpicture}[scale = 0.75]
    \begin{axis}[
      grid = both,
      legend pos=north east
    ]
    \addplot[color = red, thick] table[x = t, y = ee1] {rect2.dat};
    \addplot[color = blue!90!black, thick] table[x = t, y = ee4] {rect2.dat};
    \legend{$e_1$,$e_4$};
    \end{axis}
    \end{tikzpicture}
    \hspace*{2em}
        \begin{tikzpicture}[scale = 0.75]
    \begin{axis}[
      grid = both,
      legend pos=north east
    ]
    \addplot[color = red, thick] table[x = t, y = e1] {rect2.dat};
    \addplot[color = blue!90!black, thick] table[x = t, y = e4] {rect2.dat};
    \legend{$e_1$,$e_4$};
    \end{axis}
    \end{tikzpicture}
\\[\baselineskip]
\centering
    \begin{tikzpicture}[scale = 0.75]
    \begin{axis}[
      grid = both,
      legend pos=north east
    ]
    \addplot[color = red, thick] table[x = t, y = ee1] {rect.dat};
    \addplot[color = blue!90!black, thick] table[x = t, y = ee4] {rect.dat};
    \legend{$e_1$,$e_4$};
    \end{axis}
    \end{tikzpicture}
    \hspace*{2.5em}
        \begin{tikzpicture}[scale = 0.75]
    \begin{axis}[
      grid = both,
      legend pos=north east
    ]
    \addplot[color = red, thick] table[x = t, y = e1] {rect.dat};
    \addplot[color = blue!90!black, thick] table[x = t, y = e4] {rect.dat};
    \legend{$e_1$,$e_4$};
    \end{axis}
    \end{tikzpicture}
    \caption{Numerical solutions obtained by the trapezoidal rule applied to MNA and MONA formulations. Top: $v_{scr} = \sin(\pi t)$; bottom: $v_{scr} = \sign (\sin(\pi t))$; left: MONA; right: MNA.}
    \label{fig:instability2}
\end{figure}

\bigskip 
\noindent\textbf{Example 2: An index-1 problem with discontinuous sources.}
As our second example, we consider a full wave rectifier circuit, see Figure~\ref{fig:rectifier}, which is one of the classical components in electric and electronic devices, usually used as an AC-DC converter. The quantities of interest here are the input voltage, i.e., the potential at node $n_1$, and the output voltage, which here corresponds to the potential at node $n_4$.  
The topology of the circuit is described by the following partial incidence matrices
$$
A_L = 
\begin{pmatrix}
1 & 0 \\ 0 & 1 \\ 0 & -1 \\ 0 & 0
\end{pmatrix}
, \quad
A_R = 
\begin{pmatrix}
0 & 0 & 0 & 0 & 0 \\ -1 & 1 & 0 & 0 & 0 \\ 0 & 0 & -1 & 1 & 0 \\ 0 & -1 & 0 & -1 & 1
\end{pmatrix}
\quad\mbox{and}\quad
A_V = 
\begin{pmatrix}
1 \\ 0 \\ 0 \\ 0
\end{pmatrix}
$$
The circuit contains four diodes, which are modeled as resistors with a nonlinear voltage-current relation $j_D = 2.5 \, (\exp(4v_D) + 1)$ corresponding to the Shockley diode model. The nonlinear conductance matrix is then defined as
\begin{align*}
G:\ &\mathbb{R}^5 \to \mathbb{R}^5 \\
& v_i \mapsto j_D(v_i) ,\quad i=1,\dots ,4,\\
& v_5 \mapsto v_5/R.
\end{align*}
For simplicity, we choose $R=1$ for the remaining resistor 
Like in \cite{waltrich2010multiport}, the magnetic coupling through the transformer is modeled by the inductance matrix 
$$
L = 
\begin{pmatrix}
27.46 & 27.57 \\ 27.57 & 27.75
\end{pmatrix}
\cdot 10^{-6}.
$$
The conventianal MNA formulation here leads to an index-1 system which reads
$$
\begin{pmatrix}
0 & 0 & 0 \\ 0 & L & 0 \\ 0 & 0 & 0
\end{pmatrix}
\ddt
\begin{pmatrix}
e \\i_C \\i_V
\end{pmatrix}
+
\begin{pmatrix}
{A_R} G ({A_R}^Te){A_R}^T & A_L & A_V \\ -A_L^T & 0 & 0 \\ -A_V^T & 0 & 0 
\end{pmatrix}
\begin{pmatrix}
e \\ i_L \\ i_V
\end{pmatrix}
=
\begin{pmatrix}
0 \\0\\-v_{src}(t)
\end{pmatrix}
$$
and the MONA formulation yields a corresponding index-0 problem given by
$$
\begin{pmatrix}
{A_R} G({A_R}^T\ddt\psi){A_R}^T & A_V \\ -A_V^T & 0
\end{pmatrix}
\ddt
\begin{pmatrix}
\psi \\ q_V
\end{pmatrix}
+
\begin{pmatrix}
A_L L^{-1} A_L^T & 0 \\ 0 & 0
\end{pmatrix}
\begin{pmatrix}
\psi \\ q_V
\end{pmatrix}
=
\begin{pmatrix}
0\\-v_{src}(t)
\end{pmatrix}.
$$
Like in the previous test case, we use the trapezoidal rule for our computations with a fixed time step $\tau = 0.1$. 
In the top row of Figure~\ref{fig:instability2}, we display the interesting node potentials for a smooth input voltage $v_{src}(t) = \sin(\pi t)$ which is again consistent with trivial initial conditions. 
Both formulations are perfectly stable and yield second order convergence. 
In the second row of Figure~\ref{fig:instability2}, we display the corresponding results for the non-smooth digital input $v_{src} = \sign(\sin(\pi t))$. 
Similar to an inconsistent initial condition, the discontinuities in the source term lead to a weak instability in the MNA formulation. The MONA system, on the other hand, which here is an ordinary differential equation, is not affected by the discontinuities in the source term.

\section{Conclusion and Outlook} 

In this paper we present a magnetic oriented formulation for electrical circuits. The formulation is based on magnetic node potentials and charges across capacitors as primary unknowns, which represents an essential difference to electric-oriented approaches like the conventional MNA and charge-flux based variants. The resulting system has a generalized gradient structure which immediately ensures passivity of the system. In the majority of the cases the formulation leads to a lower index then the standard MNA formulations.
Let us remark
that the magnetic oriented viewpoint can be applied to derive a corresponding \emph{magnetic oriented loop analysis}. 
Similar to the port-Hamiltonian formulation of the MNA, the particular structure of MONA may be beneficial for multiphysical coupling required, e.g., in the context of electric machine simulation.

{\footnotesize
\section*{Acknowledgement}
Support by the German Science Foundation (DFG) via grants TRR~146 (project~C3), TRR~154 (project~C04) and SPP~2256 (project Eg-331/2-1)
and by the Center for Computational Engineering at TU~Darmstadt is gratefully acknowledged. 
}


\begin{thebibliography}{10}

\bibitem{amrhein20093}
M.~Amrhein and P.~T. Krein.
\newblock 3-d magnetic equivalent circuit framework for modeling
  electromechanical devices.
\newblock {\em IEEE Trans. Energy Conversion}, 24:397--405, 2009.

\bibitem{bannenberg2021coupling}
M.~Bannenberg, A.~Ciccazzo, and M.~G{\"u}nther.
\newblock Coupling of model order reduction and multirate techniques for
  coupled dynamical systems.
\newblock {\em Applied Mathematics Letters}, 112:106780, 2021.

\bibitem{bartel2011structural}
A.~Bartel, S.~Baumanns, and S.~Sch{\"o}ps.
\newblock Structural analysis of electrical circuits including
  magnetoquasistatic devices.
\newblock {\em Appl. Numer. Math.}, 61:1257--1270, 2011.

\bibitem{bartel2018pdaes}
A.~Bartel and M.~G{\"u}nther.
\newblock {PDAE}s in refined electrical network modeling.
\newblock {\em SIAM Rev.}, 60:56--91, 2018.

\bibitem{brenan1995numerical}
K.~E. Brenan, S.~L. Campbell, and L.~R. Petzold.
\newblock {\em Numerical solution of initial-value problems in
  differential-algebraic equations}.
\newblock SIAM, 1996.

\bibitem{brezzi1974existence}
F.~Brezzi.
\newblock On the existence, uniqueness and approximation of saddle-point
  problems arising from {L}agrangian multipliers.
\newblock {\em RAIRO}, 8:129--151, 1974.

\bibitem{campbell2019applications}
S.~Campbell, A.~Ilchmann, V.~Mehrmann, and T.~Reis.
\newblock {\em Applications of differential-algebraic equations: examples and
  benchmarks}.
\newblock Springer, 2019.

\bibitem{GarciaEtAl2020}
I.~Cortes~Garcia, H.~De~Gersem, and S.~Sch\"ops.
\newblock A structural analysis of field/circuit coupled problems based on a
  generalised circuit element.
\newblock {\em Numer. Algor.}, 83:373--394, 2020.

\bibitem{egger2019structure}
H.~Egger.
\newblock Structure preserving approximation of dissipative evolution problems.
\newblock {\em Numerische Mathematik}, 143:85--106, 2019.

\bibitem{egger2020energy}
H.~Egger, O.~Habrich, and V.~Shashkov.
\newblock On the energy stable approximation of hamiltonian and gradient
  systems.
\newblock {\em Computational Methods in Applied Mathematics}, -1, 12 2020.

\bibitem{estevez2000structural}
D.~Est{\'e}vez~Schwarz and C.~Tischendorf.
\newblock Structural analysis of electric circuits and consequences for mna.
\newblock {\em International Journal of Circuit Theory and Applications},
  28:131--162, 2000.

\bibitem{griepentrog1986}
E.~Griepentrog and R.~M\"arz.
\newblock {\em Differential-algebraic equations and their numerical treatment}.
\newblock Teubner, Leipzig, 1986.

\bibitem{gunther2021dynamic}
M.~G{\"u}nther, A.~Bartel, B.~Jacob, and T.~Reis.
\newblock Dynamic iteration schemes and port-hamiltonian formulation in coupled
  differential-algebraic equation circuit simulation.
\newblock {\em Int. J. Circuit Theory Appl.}, 49:430--452, 2021.

\bibitem{gunther1995dae}
M.~G{\"u}nther and U.~Feldmann.
\newblock The dae-index in electric circuit simulation.
\newblock {\em Math. Comput. Sim.}, 39:573--582, 1995.

\bibitem{gunther2005modelling}
M.~G{\"u}nther, U.~Feldmann, and J.~ter Maten.
\newblock Modelling and discretization of circuit problems.
\newblock {\em Handbook of numerical analysis}, 13:523--659, 2005.

\bibitem{hairer2014}
E.~Hairer and C.~Lubich.
\newblock Energy-diminishing integration of gradient systems.
\newblock {\em IMA J. Numer. Anal.}, 34:452--461, 2014.

\bibitem{hairer2006numerical}
E.~Hairer, C.~Lubich, and M.~Roche.
\newblock {\em The numerical solution of differential-algebraic systems by
  Runge-Kutta methods}, volume 1409.
\newblock Springer, 2006.

\bibitem{ho1975modified}
C.-W. Ho, A.~Ruehli, and P.~Brennan.
\newblock The modified nodal approach to network analysis.
\newblock {\em IEEE Trans. Circuits and Systems}, 22:504--509, 1975.

\bibitem{KunkelMehrmann2006}
P.~Kunkel and V.~Mehrmann.
\newblock {\em Differential-algebraic equations}.
\newblock European Mathematical Society (EMS), Z\"{u}rich, 2006.
\newblock Analysis and numerical solution.

\bibitem{laithwaite1967magnetic}
E.~R. Laithwaite.
\newblock Magnetic equivalent circuits for electrical machines.
\newblock In {\em Proceedings of the Institution of Electrical Engineers},
  volume 114, pages 1805--1809, 1967.

\bibitem{lamour2013}
R.~Lamour, R.~M\"arz, and C.~Tischendorf.
\newblock {\em Differential-algebraic equations: A Projector Based Analysis}.
\newblock Springer, Berlin, 2013.

\bibitem{petzold1982differential}
L.~Petzold.
\newblock Differential/algebraic equations are not ode’s.
\newblock {\em SIAM J. Sci. Stat. Comput.}, 3:367--384, 1982.

\bibitem{reis2014mathematical}
T.~Reis.
\newblock Mathematical modeling and analysis of nonlinear time-invariant rlc
  circuits.
\newblock In {\em Large-scale networks in engineering and life sciences}, pages
  125--198. Springer, 2014.

\bibitem{riaza2008differential}
R.~Riaza.
\newblock Differential-algebraic systems: Analytical aspects and circuit
  applications, 2008.

\bibitem{ruehli1974equivalent}
A.~E. Ruehli.
\newblock Equivalent circuit models for three-dimensional multiconductor
  systems.
\newblock {\em IEEE Trans. Microwave Theory and Techniques}, 22:216--221, 1974.

\bibitem{scholz2018signature}
L.~Scholz.
\newblock The signature method for daes arising in the modeling of electrical
  circuits.
\newblock {\em Journal of Computational and Applied Mathematics}, 332:107--139,
  2018.

\bibitem{tischendorf1999topological}
C.~Tischendorf.
\newblock Topological index calculation of differential-algebraic equations in
  circuit simulation.
\newblock {\em Surveys Math. Indust.}, 8:187--199, 1999.

\bibitem{tischendorf2003coupled}
C.~Tischendorf.
\newblock Coupled systems of differential algebraic and partial differential
  equations in circuit and device simulation.
\newblock {\em Modeling and numerical analysis}, 2003.

\bibitem{waltrich2010multiport}
G.~Waltrich, J.~Duarte, and M.~Hendrix.
\newblock Multiport converters for fast chargers of electrical vehicles-focus
  on high-frequency coaxial transformers.
\newblock In {\em The 2010 International Power Electronics Conference-ECCE
  ASIA-}, pages 3151--3157. IEEE, 2010.

\end{thebibliography}

\end{document}